\documentclass[10pt]{article}
\usepackage{a4wide}
\usepackage{wasysym}
\usepackage[english]{babel}
\usepackage[all]{xy}
\usepackage{amsfonts}
\usepackage{amssymb}
\newtheorem{defi}{\textbf{Definition}}[section]
\newtheorem{theo}[defi]{\textbf{Theorem}}
\newtheorem{lemma}[defi]{\textbf{Lemma}}
\newtheorem{prop}[defi]{\textbf{Proposition}}
\newtheorem{coro}[defi]{\textbf{Corollary}}

\newcommand{\G}{\Gamma}

\title{Parabolic groups acting on one-dimensional compact spaces}

\author{Fran\c{c}ois Dahmani \footnote{The author acknowledges support of the FIM, ETH Z\"urich.}} 

\date{}
\begin{document}

\maketitle
\begin{center}

  {\footnotesize 
{\bf Abstract.} 
Given a class of compact spaces, 
we ask which groups can be maximal parabolic subgroups of  
a relatively hyperbolic group whose boundary is in the class. 
We investigate the class of 1-dimensional connected boundaries. 
We get that any non-torsion infinite finitely generated group is a maximal parabolic 
subgroup of some relatively 
hyperbolic group with connected one-dimensional boundary without global cut point.  
For boundaries homeomorphic to a Sierpinski carpet or a 2-sphere, 
the only maximal parabolic subgroups allowed are 
virtual surface groups (hyperbolic, or virtually $\mathbb{Z} + \mathbb{Z}$). 
}
\end{center}

Let $M$ be a hyperbolic manifold of finite volume. The maximal parabolic
subgroups of $\pi_1(M)$, for its action on the universal cover of $M$, are
virtually Abelian, by the Bieberbach Theorem. 
 More generally, the parabolic
 subgroups of a geometrically finite group on a Hadamard manifold with pinched
 negative curvature are virtually nilpotent: according to \cite{Bparabolic},
 they are finitely generated, and the Margulis Lemma (Theorem 9.5 in \cite{BGS}) applies.

 A natural generalization of the class of geometrically finite groups is the class of relatively hyperbolic groups. They were first introduced by M.~Gromov \cite{G}, and studied by B.~Bowditch \cite{Brel}, and independently by B.~Farb \cite{F}.  We will follow Bowditch's approach  (see \cite{Brel} \cite{Szcz} and \cite{Dproc}(appendix) for equivalence of definitions). 
A finitely generated  group $\Gamma$ is hyperbolic relative to a family of finitely generated subgroups $\mathcal{G}$, in the sense of \cite{Brel}, if it acts properly discontinuously by isometries on a proper hyperbolic length space, such that the action induced on the boundary is a geometrically finite convergence action, whose maximal parabolic subgroups are the elements of $\mathcal{G}$. The definitions are developed in the first section below. The boundary of such a space is shown in \cite{Brel} to be canonically associated to the pair $(\Gamma, \mathcal{G})$. We call it the Bowditch boundary of the relatively hyperbolic group.

In this paper, we will see that the consequence of the Margulis Lemma mentionned above is, to a large extend, false for relatively hyperbolic groups, but may 
remain true for certain specific classes of boundaries.

In fact, it is easy to construct counterexamples, that is, relatively hyperbolic groups with an arbitrary  finitely generated parabolic subgroup $H$: it suffices to consider the 
 free HNN extension $H*_{\{1\}}$. It is hyperbolic relative to the conjugates of $H$ (see Definition 2 in \cite{Brel}, where the graph involved is the Bass-Serre tree).  However, the Bowditch boundary of this 
group is not connected. It is a Cantor set (\cite{Brel}, \cite{Dcomb}).

Our first theorem is a non-trivial generalization of this example.

\begin{theo} \label{theo;intro1}

Let $H$ be an infinite finitely generated  group which is not a torsion group. Then, there exists a relatively hyperbolic group $\Gamma$, containing $H$ 
as a maximal parabolic subgroup, and whose boundary is connected, locally connected, 1-dimensional compact space without global cut point. 
\end{theo}

The definition of relative hyperbolicity easily implies that the maximal
parabolic subgroups are infinite when the boundary has no isolated point. We
require that they are finitely generated in order to ensure the equivalence
between the different common definitions of relative hyperbolicity. We also
use this assumption explicitly to construct a suitable group $\Gamma$ whose boundary has no global cut point. It is unclear whether the assumption of not being a torsion group is necessary.

To prove the Theorem, we will make use of the 
Combination Theorem of \cite{Dcomb} for 
amalgams over infinite cyclic groups. This is done in Section 2. 
In the beginning of this section, we also give a very short proof, 
but still using the Combination Theorem, 
of the same result without the requirement of absence of global cut points. 
 In both cases, the amalgams over cyclic groups introduce local cut points in the boundary. 

In \cite{KK}, M.~Kapovich and B.~Kleiner proved that, if the boundary of a hyperbolic group is one-dimensional, connected and has no local cut point, then it is a Sierpinski carpet or a Menger curve. Their argument, recalled in section 3, remains valid  for boundaries of relatively hyperbolic groups. This gives motivation for the generalization of another of their results (Theorem 5 in \cite{KK}).

 \begin{theo}\label{theo;extension}

 Let $G$ be a relatively hyperbolic group whose Bowditch boundary is a Sierpinski carpet. 
There exists a relatively hyperbolic group $\Gamma$ whose 
boundary is a 2-sphere, and in which  $G$ embeds as a fully quasi-convex subgroup 
(in the sense of \cite{Dcomb}). 
\end{theo}

With this theorem, we will get a positive result for our original question for relatively hyperbolic groups whose boundary is a Sierpinski carpet.

\begin{theo} \label{theo;intro2}

 If  $(\Gamma,\mathcal{G})$ is a relatively hyperbolic group whose boundary is a Sierpinski carpet or a 2-sphere, 
then the maximal parabolic 
subgroups of $\Gamma$ are  virtual surface groups.  
 
\end{theo}

If only minimal families of parabolic subgroups are to be considered, 
the theorem becomes:

\begin{coro} (Corollary \ref{coro;sierpinski_minimal})

Let $(\Gamma,\mathcal{G})$ be a  relatively hyperbolic group, such that $\mathcal{G}$ is minimal for this property, and  whose boundary is a Sierpinski carpet or a 2-sphere. Then every  element of $\mathcal{G}$ is virtually $\mathbb{Z} +\mathbb{Z}$.

\end{coro}

We give, at the end of the paper, examples of relatively hyperbolic groups
with boundary homeomorphic to the Sierpinski carpet, and with parabolic
subgroups isomorphic to surface groups. 
As we mentioned, according to the result of M.~Kapovich and B.~Kleiner, if the boundary of a relatively hyperbolic group is connected one-dimensional and without local cut point, then it is homeomorphic either to the Sierpinski carpet or to the Menger curve. Thus, Theorem \ref{theo;intro1} and Theorem \ref{theo;intro2} give information on the possible parabolic subgroups of relatively hyperbolic groups whose boundaries are one-dimensional connected compact spaces, except the Menger curve. This remaining case is still  unclear to us, although we can mention that a bounded parabolic group acting on it must have one end (Proposition \ref{prop;menger}).  I do not know an example of a one-ended finitely generated group which is not a parabolic subgroup of a relatively hyperbolic group
whose boundary is a Menger curve.

\vskip .3cm
I would like to thank A.~Szczepanski for the discussions we had about related problems, and S.~Maillot, who helped me to simplify a part of the proof. I am indebted to the referee for helpful remarks, advice, and suggestions.

\section{About relatively hyperbolic groups}

We briefly set the framework. 
The concept of relatively hyperbolic group has been suggested by M.~Gromov, 
and elaborated on by B.~Bowditch and independently by B.~Farb. 
For more details and different equivalent 
definitions, see \cite{Brel} \cite{F} \cite{G}.

We recall that a group $G$ acts on a compact space $K$ as a 
\emph{convergence group}, if it acts 
properly discontinuously 
on the space of distinct triples of $K$. 
A point $\xi \in K$ is a \emph{conical limit point} if there exists a 
sequence $g_n$ of elements 
of $G$ and two distinct elements $a \neq b$ in $K$ such that $g_n \xi \to a$ 
and $g_n \zeta \to b$ for all 
$\zeta \in K \setminus \{ \xi \}$. Finally, $\xi \in K$ is a 
\emph{bounded parabolic point} 
if its stabilizer $Stab(\xi)$ 
acts properly discontinuously and 
co-compactly on $K\setminus\{\xi\}$ (see \cite{BM} \cite{Brel} \cite{TukiaCrelle}). 
The action of a convergence group on a compact space is \emph{geometrically finite} if this compact space consists only of conical limit points and bounded parabolic points. The maximal parabolic subgroups are then the stabilizers of the parabolic points.

\begin{defi}(Relatively hyperbolic groups) \cite{Brel}

Let $\Gamma$ be a finitely generated group, and $\mathcal{G}$ be a (closed by conjugation) family of finitely generated subgroups of $\Gamma$. We say that $(\Gamma, \mathcal{G})$ is a relatively hyperbolic 
group if $\Gamma$ acts properly discontinuously by isometries on a proper Gromov-hyperbolic length space $X$, inducing an action by homeomorphisms on the boundary $\partial X$ that makes $\Gamma$ a geometrically finite convergence group on $\partial X$, 
whose maximal parabolic subgroups are precisely the elements of $\mathcal{G}$. 
\end{defi}

The compact space $\partial X$, and the action of $\Gamma$ on it, is
canonically associated to $(\Gamma, \mathcal{G})$, and $\partial X$ is called 
 the Bowditch boundary of the relatively hyperbolic group, and is denoted by
 $\partial \Gamma$ when no space $X$ is explicitly introduced. 

The first
 examples of relatively hyperbolic groups are fundamental groups of finite
 volume (or, more generally, geometrically finite) hyperbolic manifolds. One
 can also consider variable curvature: geometrically finite groups on pinched Hadamard manifolds are hyperbolic relative to their maximal parabolic subgroups (see \cite{Brel} for more details). Other examples are free
 products: if $A$ and $B$ are infinite, finitely generated groups, the group
 $A*B$ is hyperbolic relative to the conjugates of $A$ and $B$. Limit groups,
 or finitely generated fully residually free groups, are also relatively
 hyperbolic (see \cite{Dcomb}).

\begin{defi}(Fully quasi-convex subgroups) \cite{Dcomb}\label{defi;fqc}

Let $(G,\mathcal{G})$ and $(H,\mathcal{H})$ be two relatively hyperbolic groups, and assume that $H$ is a subgroup of $G$. 
It is  quasi-convex if its limit set $\Lambda H \subset \partial G$ is equivariantly homeomorphic to its boundary $\partial H$. 
It is fully quasi-convex if moreover any infinite family of $G/H$-translates of $\Lambda H$ intersect trivially together. 
\end{defi}

We will make use of a special case of the Combination Theorem of \cite{Dcomb}, and we 
recall here the statement we need.

If a group $G$ act on a tree $T$, we say, following Z.~Sela, that the action
 is \emph{acylindrical} if there exists a number $k$ such that for all segment
 of length $k$ in $T$, its stabilizer in $G$ is finite. 
 Let us also notice that any finitely generated group is hyperbolic relative to the family consisting of itself, and that its boundary is then a single point.

\begin{theo}(\cite{Dcomb})\label{theo;comb}

  {\bf 1.} Let $\Gamma$ be the fundamental group of an acylindrical finite graph
of relatively hyperbolic groups, whose edge groups are fully
quasi-convex subgroups of the adjacent vertices groups. Let $\mathcal{G}$ be the
family of the images of the maximal parabolic subgroups of the vertices groups,
and their conjugates in $\Gamma$.
 Then $(\Gamma,\mathcal{G})$ is a relatively hyperbolic group.

 {\bf 2.} Let $G$ be a group which is hyperbolic relative to a family of
subgroups $\mathcal{G}$, and let $ P$ be a group in  $\mathcal{G}$. Let $A$ be a
finitely generated group in which $P$ embeds as a subgroup. Then, $\Gamma = A*_P
G$ is hyperbolic relative to the family $(\mathcal{H} \cup
\mathcal{A})$, where   $\mathcal{H}$  is the set of the conjugates of the
images of elements of $\mathcal{G}$ not conjugated to $P$ in $G$, and
where $\mathcal{A}$ is the set of the conjugates of $A$ in $\Gamma$.

Moreover, in both cases, if the
topological  dimensions of the boundaries of the vertex groups (resp.~of the
edge groups) are smaller than $r$ (resp.~than $s$), then
$\hbox{dim}(\partial\Gamma) \leq \max \{r,s+1\}$.

Finally, if the boundaries of every vertex groups are connected, and if the boundaries of every edge groups are non empty, then the boundary of $\Gamma$ is connected. 

\end{theo} 

The last assertion is not addessed in \cite{Dcomb}. However, it is an easy corollary of the construction of the boundary in \cite{Dcomb}. Indeed, let $T$ be the Bass-Serre tree and $\partial T$ its boundary, let $\Omega$ be the disjoint union of the boundaries of the vertex stabilizers in $T$, and $\sim$ be the equivalence relation obtained by gluing together the boundaries of the stabilizers of any two adjacent vertices along the limit set of the corresponding edge stabilizer. This quotient is Hausdorff because the equivalence classes are closed subsets, and, for any distance compatible for the topology, for all $\epsilon>0$, only finitely many equivalence classes have diameter greater than $\epsilon$. 
 Then, it is shown in \cite{Dcomb} that $\partial \Gamma$ is equivariantly homeomorphic to $(\partial T \cup \Omega) /_{\sim}$. With the hypothesis of the last assertion, it is easy to see that $\Omega /_{\sim}$ is connected, and it is dense in $\partial \Gamma$.

\section{Examples of relatively hyperbolic groups with one-dimensional boundaries.}

This section is devoted to a proof of Theorem \ref{theo;intro1}.

We begin by constructing a group $\Gamma_c$ with the same properties than 
the group to contruct, without 
the requirement on global cut points (the index $c$ stands for ``cut points''). 
As this is less technical, we indicate in the process of the construction the point where the reader only interested in existence of relatively hyperbolic groups with connected boundaries with a given parabolic subgroup may stop: this is Proposition \ref{prop;sans_coupure_global}. 
This first construction will be done easily with the Combination Theorem, 
but we keep a trace of the behaviour of the cut points.
 Then we will perform a little more delicate construction in order to get a 
boundary without cut point.

\vskip .3cm

{\it Proof of Theorem \ref{theo;intro1}. } We fix a finitely generated group $H$ that has an element $h$ of infinite order.  

 Let $F$ be a free group of rank $2$, $F=F(a,b)$, and let $\mathcal{F}$ consists of the conjugates of the cyclic subgroup generated by the commutator $[a,b]$. It is  well known (see \cite{F}) that $(F,\mathcal{F})$ is a relatively hyperbolic group, and that its boundary $\partial F$ is a circle. 
We choose $\Gamma_c$ to be the amalgamated free product $H *_{\mathbb{Z}} F$, where $\mathbb{Z}$ is 
identified to $\langle h \rangle$ in $H$, and to $\langle [a,b] \rangle$ in $F$. By the Combination Theorem \ref{theo;comb}, 
it is hyperbolic relative to the family of the conjugates of $H$, its boundary is connected,  
and has topological dimension 1.  This proves the proposition:

\begin{prop}\label{prop;sans_coupure_global}
Let $H$ be an infinite  finitely generated  group which is not a 
torsion group. Then, there exists a relatively hyperbolic group $\Gamma_c$, containing $H$ 
as a maximal parabolic subgroup, and whose boundary is a connected 1-dimensional compact space. 
\end{prop}

 Let us describe  more precisely the boundary of $\Gamma_c$. 
 Let $T$ be the Bass-Serre tree of the amalgam. In order to avoid confusion with another graph of groups to come, we denote the vertices and the edges of $T$ by Greek letters ($\nu$ and $\epsilon$).  Let  $\nu_F$ and $\nu_H$ 
be two adjacent vertices in $T$  stabilized respectively by $F$ and $H$.  
 One has 
$$ \partial \Gamma_c \simeq \left( \partial T \: \cup \: \bigcup_{\nu \in \Gamma_c \nu_F} \partial (Stab(\nu)) \cup \Gamma_c \nu_H \right)/\sim $$
 where the relation $\sim$ is the equivalence relation generated by the identification, for any adjacent vertices $\gamma \nu_H$ and $\gamma \nu_F$, of the point $\gamma \nu_H$ with the parabolic point of $\partial (Stab(\gamma \nu_F))$ fixed by the stabilizer of the edge between $\gamma \nu_H$ and $\gamma \nu_F$.  
 By an abuse of notation, we will identify $\Gamma_c \nu_H$ to the set of parabolic points in $\partial \Gamma_c$, that is, to its image after the identification $\sim$ (note that this quotient is injective on the set $\Gamma_c \nu_H$).

Let $(\epsilon_i)_{i\in I}$ be the set of edges adjacent to $\nu_H$, and for every index $i\in I$, 
denote by 
$T_{i}$ the connected component of $T\setminus \{\nu_H \}$ that contains the vertex of $\epsilon_i \setminus \{ \nu_H \}$. Let $C_{i}\subset \partial \Gamma_c$ be the image of $\partial T_{i}  \cup \left(  \bigcup_{\nu\in \Gamma_c \nu_F \cap T_{i}} \partial (Stab(\nu))\right) \cup (\Gamma_c \nu_H \cap T_{i}) $ by the canonical  projection in the quotient by the relation $\sim $. 
Note that  $C_i$ contains the parabolic point $\nu_H$ since this latter is identified to 
a parabolic point of the boundary of the stabilizer of the vertex of $\epsilon_i\setminus \{  \nu_H \}$.

It is easily seen that the set of global cut points of $ \partial \Gamma_c$ is the set 
of the
parabolic points, and that given a parabolic point $\gamma \nu_H$, the connected components of $\partial \Gamma_c \setminus \{ \gamma \nu_H \}$ are precisely the sets $\gamma (C_i \setminus \{\nu_H\})$ for $i \in I$.

\begin{lemma} \label{lem;loxo}

 Let $\epsilon_i$ and $\epsilon_j$ be two distinct edges of $T$ adjacent to $\nu_H$. Let $h\in H$ such that $h \epsilon_i = \epsilon_j$, and let $\sigma$ be an element of $\Gamma_c$ not in $H$, fixing the other vertex of $\epsilon_j$. Then $\gamma = \sigma h$ is a hyperbolic isometry of $T$, it has exactly two fixed points in  
 $\partial \Gamma_c$, one of them is in $C_i$ while the other is in $C_j$.  

\end{lemma}

{\it Proof. } Let $\nu_i$ (resp.~$\nu_j$) be the vertex of $\epsilon_i$ (resp.~$\epsilon_j$), different from $\nu_H$. Then $(\sigma h) \nu_i = \nu_j$. Therefore, $(\sigma h)$ does not fix the middle of the segment $[\nu_j, (\sigma h) \nu_i]$, which is $\nu_H$. But any elliptic isometry $i$ of a tree $T$ fixes the middle of every segment $[\nu,i(\nu)]$, $\nu\in T$. This implies that $(\sigma h)$ is not elliptic, hence it is a hyperbolic isometry. Moreover, its translation length is at most $ dist(\nu_i, 
(\sigma h)\nu_i)=2$, and it cannot be $1$, since parity is preserved. 
Therefore the bi-infinite geodesic fixed by   $(\sigma h)$ contains both  $\nu_i$ and $\nu_j$. $\square\smallskip$

 We now define another group. Let $\{h_1 \dots h_s\}$ be a symetric generating set of $H$. Let us choose $\sigma$, an element of $\Gamma_c$, not in $H$, that fixes a neighbour of $\nu_H$ in $T$.

Let $X$ be the graph of groups consisting 
in $s+1$ vertices $V_0, \dots, V_s$, of $s$ edges $e_k =(V_0,V_k)$, $k=1,\dots, s$, as follows. 
The group of $V_0$ is $\Gamma_c$, the group of every other vertex $V_k$ is a genus $2$ closed surface group $\Sigma_k$, and, for every $k$,  
the group of $e_k$ embeds as an infinite cyclic, maximal cyclic subgroup in $\Sigma_k$, and as  $\langle \sigma h_i \rangle$ in the group $\Gamma_c$.
 By the Combination Theorem, the group $\pi_1(X)$ (defined up to conjugacy) is hyperbolic relative to the family of the images of the maximal 
parabolic groups of $\Gamma_c$, and their conjugates. Let $M$ be its boundary. It is a connected compact space 
of dimension 1 by the Combination Theorem. 
More precisely, let $\tilde{X}$ be the Bass-Serre tree of $X$, and $\partial \tilde{X}$ be its boundary. The boundary of $\pi_1(X)$ is homeomorphic to  
$$M \simeq \left( \partial \tilde{X}  \: \cup \: (\bigcup_{v \in \tilde{X}^0} \partial (Stab(v))  \right) /\approx $$
 where the relation $\approx$ is the equivalence relation generated by the identification, for any adjacent vertices $v$, and $v'$, of the limit set in $\partial (Stab(v))$ and $\partial (Stab(v))$ of the stabilizer of the edge $(v,v')$.

Let $v_{\Gamma_c}$ be the vertex of $\tilde{X}$ stabilized by $\Gamma_c$.  We choose an equivariant homeomorphism  (hence we will not distinguish in notations) between its limit set in $M$ and the compact set $\partial \Gamma_c$: $\Lambda (Stab(v_{\Gamma_c})) \simeq \partial (Stab( v_{\Gamma_c}))=\partial \Gamma_c$.

{\it Remark.}  Let us remark that $M$ is locally connected. It follows (as we explain now) from the description of the fundamental system of neighborhoods given in \cite{Dcomb}(par. 2.3). Indeed, if a point $\xi$ is the image of a point of the boundary of $\tilde{X}$, the Bass-Serre tree of $X$, then a system of neighborhoods  of  $\xi$  consists of  all the boundaries of the groups fixing vertices in a component of $\tilde{X} \setminus \{e_n\}$, where $e_n$ is an edge on a ray going to $\xi$.  Since every edge group has non-empty boundary, and every vertex group has connected boundary, this gives a connected set for all $n$.
  
  If now $\xi$ is the image of a point in the boundary of a group fixing a vertex $v$, then possibly it is in the boundary of the stabilizer of an edge $(v,v')$, or may be it is not (in such case we note $v'=v$). Then a neighborhood of $\xi$ in the fundamental system proposed in \cite{Dcomb}(par. 2.3) is $W_n(\xi)  = A_n\cup B_n\cup C_n = \overline{B_n\cup C_n}$, where  $C_n$ is defined to be the union of a neighborhood of $\xi$ in   $\partial (Stab( v))$ and in   $\partial (Stab( v'))$, in fundamental basis of neighborhoods in each. By local connectivity of the boundaries of stabilizers of vertices, one can choose $C_n$ to be connected. Finally, $B_n$ is defined to be    a union over some components of $\tilde{X} \setminus{v,v'} $,     of all the boundaries of the groups fixing vertices in that component, noted $B_n= \cup_{i\in I_n}  B_i$, where $I_n$ contains the indices of the components so that  $B_i$ intersects $C_n$. Since every edge group has non-empty boundary, and every vertex group has connected boundary, each $B_i$ is connected.  Thus $B_n \cup C_n$ is connected, and so  is  $W_n(\xi)$. 

\begin{lemma} \label{lem;xi}

Let $\xi$ in $M$ be a global cut point of $M$. Then $\xi$ is a parabolic point. 

\end{lemma}

{\it Proof. } The boundary is connected, and locally connected (see remark above). The lemma is then a consequence of a theorem of Bowditch (Theorem 9.2 in \cite{BPeriph}) asserting that if the boundary of a relatively hyperbolic group is connected and locally connected, any global cut point is fixed by a vertex group of the maximal peripheral splitting, hence is a parabolic point. $\square \smallskip$

 However, in our case, it is easy to check that any point in $\partial \tilde{X}$, and any point in the limit set of $\partial Stab(v_{\Gamma_c})$ that is not a parabolic point, is not a global cut point.

\begin{lemma} \label{lem;2.4}

 Let $\xi$ be a parabolic point of $M$ in $\partial (Stab(v_{\Gamma_c}))$.  
Let $\mathcal{O}$ be a non-empty open closed subset of
 $M\setminus\{ \xi \}$. Then $\mathcal{O}$ contains a point  of $\partial (Stab (v_{\Gamma_c}))$.

Moreover, let $v$ and $w$ be two points in $\tilde{X} $ such that the segment $[v,w]$ does not contain $v_{\Gamma_c}$.  If $\mathcal{O}$ contains a point in $\partial (Stab(v))$, 
then it contains $\partial (Stab(v))$ and  $\partial (Stab(w))$.

\end{lemma}

{\it Proof. } 
 Let us remark that, $\xi$ being a parabolic point in $\partial(Stab(v_{\Gamma_c}))$, it belongs to $\partial Stab (v)$ only if $v=v_{\Gamma_c}$.
Let us prove the second assertion first. The set $\mathcal{O}$ is open-closed, and  $\partial (Stab (v))$ is connected, therefore  $\partial (Stab (v))\subset \mathcal{O}$. Let us consider a vertex $w$ such that the segment $[v,w]$ does not contain $v_{\Gamma_c}$ in $\tilde{X}$. We denote by $e_1, \dots, e_l$ the consecutive edges of this segment, and its consecutive vertices are $v,v_1,\dots, v_l=w$.   The set $\mathcal{O}$ intersects nontrivially $\partial(Stab(v_1))$ since it intersects $\partial(Stab(e_1))$. 
Since  $\partial(Stab(v_1))$ is connected, and does not contain $\xi$ (because $v_1\neq v_{\Gamma_c}$), one has  $\partial(Stab(v_1)) \subset \mathcal{O}$. This iteration can be done $l$ times in order to get that $\mathcal{O} \cap \partial (Stab (e_l)) \neq \emptyset$, and therefore, by connectedness, $\partial (Stab(w)) \subset \mathcal{O}$.

Let us now prove the first assertion. Let $\zeta $ be a point in $\mathcal{O}$, and assume that $\zeta \notin \partial (Stab(v_{\Gamma_c}))$. There are two possibilities. Either $\zeta \in \partial \tilde{X}$, or $\zeta \in \partial (Stab(v))$ where $v$ is a vertex of $\tilde{X}$ different from $ v_{\Gamma_c} $. 

If we are in the first case, as $\mathcal{O}$ is open, we can deduce that it contains the limit set of the stabilizers of the vertices that are in some neighbourhood of $\zeta$. Therefore, we can assume without loss of generality that we are in the second case. 

From the second assertion, we deduce that there exists $v$ a neighbour of $v_{\Gamma_c}$ such that $\partial (Stab(v)) \subset \mathcal{O}$. As  $\partial (Stab(v)) \cap \partial (Stab(v_{\Gamma_c}))$ consists of the two points that are fixed by the edge between them (these points are not parabolic by choice of the edge groups),  $\mathcal{O} \cap \partial (Stab (v_{\Gamma_c}))$ contains 
these two points. $\square\smallskip$

\begin{prop} \label{prop;pas_xi}

  Let $\xi$ be a parabolic point of $M$. Then $M \setminus  \{\xi\}$ is connected. 
\end{prop}

{ \it  Proof. } After composing by an element of $\pi_1(X)$, one can assume
that $\xi \in \partial(Stab(v_{\Gamma_c}))$.  By the previous lemma, it is enough to show that if
$\mathcal{O}$ is open-closed in $M\setminus \{\xi\}$, and   contains a connected component of $\partial (Stab (v_{\Gamma_c})) \setminus \{\xi \}$, then it contains them all. Indeed if it is the case, by Lemma \ref{lem;2.4}, there is only one such set $\mathcal{O}$ possible, hence it has to be $\mathcal{O} = M \setminus \{\xi\}$.

Therefore it is enough to show the next lemma. 

\begin{lemma}

 Let $C$  and $C'$ be  two connected components of $\partial (Stab (v_{\Gamma_c})) \setminus \{\xi \}$. 
There exists an integer $n$, a sequence  $e_1, \dots, e_n$ of edges of $\tilde{X}$ adjacent to $v_{\Gamma_c}$, and a sequence $C= U_0, U_1, \dots, U_n =C'$ of connected components of  $\partial (Stab (v_{\Gamma_c})) \setminus \{\xi \}$, such that for all $r=1\dots n$, $\partial (Stab(e_r)) \cap U_{r-1} \neq \emptyset$ and $\partial (Stab(e_r)) \cap U_r \neq \emptyset$.

\end{lemma}

{\it  Proof. } We can identify (by a homeomorphism) the space $\partial (Stab
 (v_{\Gamma_c})) \setminus \{\xi \}$ to $\partial \Gamma_c \setminus
 \{\nu_H\}$. Hence its connected components are the $(C_i)_{i\in I}$, indexed by the edges of $T$ (the Bass-Serre tree defining $\Gamma_c$, we do \emph{not} mean $\tilde{X}$) adjacent to the vertex $\nu_H$. 
 Let $\epsilon_{\sigma}$ be the edge adjacent to $\nu_H$ whose other end is fixed by $\sigma$, the element used in the definition of the graph $X$. By transitivity, it is sufficient to prove the result for $C'$ being the component associated to $\epsilon_{\sigma}$.

 Let $\epsilon$ be the edge of $T$ associated to $C$, and let $h\in H$ be such that $h \epsilon_{\sigma} = \epsilon$. We write $h=h_{(1)} \dots h_{(n)}$, where each $h_{(r)}$, is an element of the symmetric generating family $h_1 \dots h_s$ of $H$ used in the definition of the graph $X$.
 We now choose $\epsilon_i =  h_{(i)} h_{(i+1)} \dots h_{(n)} \epsilon$, and  $U_i$ to be $C_i \setminus \{\nu_H\}$, that is the component of $\partial \Gamma_c \setminus \{v_H\}$ associated to $\epsilon_i$. Hence we have $h_{(i)} \epsilon_{i+1} = \epsilon_i$ for all $i\leq n$. For all $i$, let $\rho_i =   h_{(i)} h_{(i+1)} \dots h_{(n)}$. Note that $\rho_i \sigma \rho_i^{-1}$ fixes the vertex of $\epsilon_i$ different from $\nu_H$. Note also that  $\rho_i h_{(i)} \rho_i^{-1} (\rho_i\epsilon_{i+1}) = \rho_i \epsilon_i$. Therefore, by Lemma \ref{lem;loxo} the element $\rho_i\sigma h_{(i)} \rho_i^{-1}$ is hyperbolic in $T$ and has one fixed point in $U_i$ and another in $U_{i+1}$. Moreover, by choice of the edge groups of $X$, the element $\rho_i\sigma h_{(i)} \rho_i^{-1}$ fixes an edge of $\tilde{X}$, since it is conjugated to  $\sigma h_{(i)}$. This proves the Lemma.  $\square\smallskip$

\begin{coro} \label{coro;nocutpoint}
 The space $M$ has no global cut point.
\end{coro}

{\it Proof. } It is a consequence of Lemma \ref{lem;xi} and Proposition
\ref{prop;pas_xi}. $\square \smallskip$   
\vskip  .5cm

{\it End of the proof of Theorem \ref{theo;intro1}. } The group $\pi_1(X)$ is
relatively hyperbolic, contains a maximal parabolic subgroup isomorphic to $H$,
and its boundary is connected, and one dimensional,  by the Combination Theorem.  By Corollary \ref{coro;nocutpoint}, its boundary has no
global cut point. 
 We already noticed that the boundary is locally connected (see remark before Lemma 2.3). $\square$

\section{Boundaries without local cut point, and parabolic groups acting on Sierpinski carpets.}

\subsection{Remarks about boundaries without local cut point}

Theorem \ref{theo;intro1} involved free constructions 
and amalgams that introduce 
local cut points in the boundary. 
We address a few remarks on what can be said if one forbid these local cut points. 
We note first that the class of boundaries considered consists in fact only of the Sierpinski carpet, 
and the Menger curve, and we give a restriction for the maximal parabolic subgroups.

Let us recall that the Sierpinski carpet is obtained as follows. Start from
 the unit square of the Euclidean plane, and divide it into nine smaller
 squares, and remove the interior of the central one. Perform this operation recursively on each of the remaining eight squares. The resulting compact space is the Sierpinski carpet. A Sierpinski curve is a compact space homeomorphic to the Sierpinski carpet.
 Let us recall the following result: 

\begin{theo} \cite{ref45deKK}\label{theo;caract_Sier}
 A compact space $\Sigma$ is homeomorphic to the Sierpinski carpet if and only if it is 1-dimensional, planar, connected and locally connected, with no local cut point.

Let $\Sigma$ be the Sierpinski carpet, then the non separating embedded circles (called ``peripheral circles'') are pairwise disjoint, and form an infinite countable set.

Given any metric (compatible with the topology) on $\Sigma$, and any number $\epsilon >0$, there are only finitely many peripheral circles of diameter greater than $\epsilon$.
\end{theo}

We recall also a result of M.~Kapovich and B.~Kleiner (stated originally for hyperbolic groups).

\begin{prop}(Theorem 4 in \cite{KK}) \label{prop;fromKK}

If the boundary of a relatively hyperbolic group is 1-dimensional, connected, and has no local cut point, then it is either a Sierpinski curve, or a Menger curve.
\end{prop}

{\it Proof. } We briefly reproduce the proof of \cite{KK} slightly modified for our context. 
First, from \cite{BGF} (Corollary 0.2), if the boundary of a relatively hyperbolic group is connected, and has no global cut point, it is locally connected.  
A compact metrisable connected, locally connected 1-dimensional space is a Menger curve 
provided it has no local cut point, and no non-empty open subset is planar. 
If the boundary $\partial \Gamma$ of a relatively hyperbolic group has a planar non-empty open 
subset $\mathcal{O}$, then any subset $S$ of $\partial \Gamma$ homeomorphic to a finite graph is in fact a planar graph. To see this, note that  $\mathcal{O}\setminus S$ is open and non-empty, since $\partial \Gamma$ has no local cut point. As the set of parabolic points is dense in $\partial \Gamma$, the open set  $\mathcal{O}\setminus S$ contains  a parabolic point $\xi$. 
 Its stabilizer acts co-compactly on $\partial \Gamma \setminus \{\xi\}$, 
therefore, there exists an element $\gamma$ such that $\gamma S \subset \mathcal{O}$, and 
therefore $S$ is planar.
 Moreover, any compact, metrisable, connected, locally connected space without cut point, and 
without non-planar embedded graph, is planar. Therefore $\partial \Gamma$ is planar, 
and by the characterisation of 
Theorem \ref{theo;caract_Sier}, 
it is homeomorphic to the Sierpinski carpet. $\square$

\begin{prop} \label{prop;menger}
 If $\Gamma$ is a relatively hyperbolic group whose boundary is a Menger curve $M$, and if $H$ is a maximal parabolic subgroup of $\Gamma$, then $H$ is one-ended.
\end{prop}

{\it Proof .} By definition, $H$ acts properly discontinuously and co-compactly on 
$M \setminus \{\xi \}$, for some point $\xi$. Since $M$ has no global cut point, this space is connected, locally connected. 
 Moreover, $M$ has no local cut point. In other words, there is a fundamental system of neighbourhoods of $\xi$, say $U_n$ 
such that each $U_n\setminus \{\xi\}$  is connected. 
Therefore, $M \setminus \{\xi\}$ is connected at infinity, in other words it has only one end.
 The group $H$ acts properly discontinuously co-compactly on it. Let $D$ be a compact fundamental domain for the action, and $N$ be the nerve of the covering of  $M \setminus \{\xi\}$ by the translates by $H$ of $D$. The group $H$ with any word metric is quasi-isometric to the complex $N$ (with the natural metric of simplicial complex). Moreover, if some compact sub-complex $C$ of $N$ separates $N$ into several unbounded component, then the union $\bigcup_{h\in H, \{h D\}\in C} h D$, of the translates of $D$ in $M \setminus \{\xi\}$ that belongs to $C$, separates   $M \setminus \{\xi\}$ into several unbounded components.  Therefore, 
$N$ has only one end, and $H$ also, because the number of ends is invariant by quasi-isometries for length spaces. $\square \smallskip$

We now prove Theorem \ref{theo;intro2}. First we need Theorem \ref{theo;extension}.

\subsection{Doubling Sierpinski carpets}

In this part, we prove Theorem \ref{theo;extension}. We follow some points of \cite{KK}.

We assume now that a group $\Gamma$ is hyperbolic relative to a family of subgroups $\mathcal{G}$, and that $\partial \Gamma$ is homeomorphic to the Sierpinski 
carpet: $\partial \Gamma \simeq \Sigma$. 
We refer to the paper of M.~Kapovich and B.~Kleiner \cite{KK} for a similar study 
for hyperbolic groups.

\begin{prop}
 
There are only finitely many $\Gamma$-orbits of peripheral circles in $\Sigma$. 
\end{prop}

{\it Proof. } In \cite{KK}, the authors use the co-compact action of the space of distinct triples; 
we could not adapt this, but a remark on the expansivity of the action gives another proof.  
We fix a metric on $\Sigma$ that is compatible with the topology. In \cite{Dsym} (Proposition 3.18), we proved that 
the action of $\Gamma$ is expansive on $\Sigma$: there exists $\epsilon >0$ 
such that for any pair of 
distinct points $x,x'$ of $\Sigma$, there exists $\gamma \in \Gamma$ such that $\gamma x$ and $\gamma x'$ are at distance at least $\epsilon$.    Given $\epsilon >0$, there are only 
finitely many peripheral circles of diameter greater than $\epsilon$. Let us note as well that every 
element $\gamma\in \Gamma $ sends a peripheral circle on a peripheral circle, since they are precisely the 
non-separating topological circles of $\Sigma$, a property that must be preserved by homeomorphisms. 
By definition of expansivity, there exists $\epsilon >0$ such that any two points 
(in particular, on a same peripheral circle $C$) 
can be cast $\epsilon$-away from each other by an element $\gamma \in \Gamma$. 
This element $\gamma$ sends $C$ onto one of the  peripheral circles of diameter greater than $\epsilon$. $\square \smallskip$

\begin{lemma}\label{lem;conical_in_C}
 
 Let $C$ be a peripheral circle of $\Sigma$, and $G < \Gamma$ be its stabilizer. Then $G$ acts as a convergence group on $C$. Moreover, any point in $C$ that is a conical limit point for $\Gamma$ is a conical limit point for $G$.   
\end{lemma}

{\it Proof. } The group $\Gamma$ acts as a convergence group on $\Sigma$, therefore, the subgroup $G$ of $\Gamma$ is of convergence on $C$. Without loss of generality (up to conjugacy), we can choose $C$ to be of maximal diameter in its orbit under the $\Gamma$-action. Let $\xi$ be a conical limit point for $\Gamma$ in $C$. 
Then there exists a sequence of elements $\gamma_n$ in $\Gamma$, and two different points $a$ and $b$ in $\Sigma$ 
such that $\gamma_n \xi \to a$ and $\gamma_n \zeta \to b$ for all  $\zeta \neq \xi$. Now note that $\gamma_n C$ ranges over only finitely many peripheral circles: if not, there would be a subsequence such that the diameter of $\gamma_{n_k} C$ collapses to zero, and then for any point $\zeta$ in $C$, $\gamma_{n_k} \zeta$ would have same limit as $\gamma_{n_k} \xi$, which is not permitted. Therefore, there is a subsequence such that $\gamma_{n_k} C$ is the same peripheral circle for all $n_k$. Because it is closed, the points $a$ and $b$ are in this circle. We translate by $\gamma_{n_0}^{-1}$, so that   $\gamma_{n_0}^{-1}  \gamma_{n_k} C = C$ for all $k$. The sequence $\gamma_{n_0}^{-1}  \gamma_{n_k} \xi$  converges to $\gamma_{n_0}^{-1} a \in C$ and for all other $\zeta \in C$,  $\gamma_{n_0}^{-1}  \gamma_{n_k} \zeta$ 
converges to  $\gamma_{n_0}^{-1} b \in C\setminus \{a\}$. For all $k$,  $\gamma_{n_0}^{-1}  \gamma_{n_k}$ is an element of $G$, because it sends $C$ on itself. This proves that $\xi$ is a conical limit point of $C$ for the action of $G$. $\square \smallskip$

\begin{lemma}\label{lem;parab_in_C}

Let $C$ be a peripheral circle of $\Sigma$, and $G<\Gamma$ be its stabilizer. 
Then $G$ is a virtual surface group, and it acts on $C$ as a uniform convergence group.
\end{lemma}

{\it Proof. } From \cite{Gab}, it is enough to prove the second assertion, and from 
\cite{TukiaCrelle}, it is enough to prove that $C$ consists only of conical limit points. 
From Lemma \ref{lem;conical_in_C}, it suffices to show that $\Sigma$ has no bounded parabolic 
point lying on $C$, for the action of $\Gamma$. Assume the contrary:~let $\xi \in C$ 
be a bounded parabolic point for $\Gamma$. Let $Stab(\xi) <\Gamma$ be its stabilizer. 
Like in the proof of Proposition \ref{prop;menger}, since $\Sigma$ is connected, locally connected, and has no cut point,  
$\Sigma\setminus \{\xi\}$ is 
connected, and connected at infinity, and therefore $Stab(\xi)$ has only one end.

Now note that $Stab(\xi)$ stabilizes the circle $C$, because two distinct peripheral circles are disjoint. Therefore, it acts 
properly discontinuously and co-compactly on $C\setminus \{\xi\}$, which is homeomorphic to the real line (hence it has two ends). 
This is a contradiction. $\square \smallskip$

\begin{prop}

The stabilizer of a peripheral circle is a fully quasi-convex subgroup of $\Gamma$ (see Def. \ref{defi;fqc}).

\end{prop}

{\it Proof. } Let $G$ be the stabilizer of a peripheral circle $C$. The Lemmas \ref{lem;conical_in_C} 
and \ref{lem;parab_in_C} show that $G$ acts as a uniform convergence group on $C$, which is clearly its limit set. 
This shows that $G$ is quasi-convex in $\Gamma$ (see Definition \ref{defi;fqc}). 
It is fully quasi-convex because  the limit set of a 
conjugate of $C$ is also a peripheral circle, and two distinct peripheral circles have empty intersection. 
$\square \smallskip$ 

We can now prove a result which is analogous to a theorem of Kapovich and Kleiner \cite{KK}.

\begin{theo}

Let $\Gamma,\mathcal{G}$ be a relatively hyperbolic group, whose boundary is homeomorphic to the 
Sierpinski Carpet: $\partial \Gamma \simeq \Sigma$. Let $C_1, \dots, C_k$ be a set of representatives of the set of peripheral 
circles of $\partial \Gamma$ under the action of $\Gamma$, and $H_1, \dots, H_k$ be respectively their stabilizers in $\Gamma$.

Let $X$ be the graph of groups consisting of two vertices $v_1$ and $v_2$, and $k$ edges $e_1, \dots, e_k$ each of 
them being between $v_1$ and $v_2$. 
The group of the vertex $v_i$ is $\Gamma$, and the group of the edge $e_i$ is $H_i$, the maps being the inclusion maps.

Then, the group $\pi_1(X)$ is hyperbolic relative to the conjugates of the images of the parabolic groups of $\Gamma$ in both sides, 
and its boundary is homeomorphic to a 2-sphere. 
\end{theo}  

The conjugate of a subgroup $H_i$ of $\Gamma$ by an element that is not in $\Gamma$ has finite intersection with $H_i$, since this subgroup has an empty limit set. This implies that the graph of groups $X$ is acylindrical. Therefore, the Combination Theorem \ref{theo;comb} gives the relative hyperbolicity.   
The boundary of $\pi_1(X)$ is completely described (see \cite{Dcomb}) by the Bass-Serre tree, 
and the boundaries of each of the vertex groups: 
$\partial \pi_1(X) \simeq (\partial T \bigcup \Omega/_\sim )$, 
where $T$ is the Bass-Serre tree, $\Omega$ is the disjoint union of the boundaries of the 
vertices stabilizers of $T$, 
 and $\sim$ is the equivalence relation defined by the attaching maps between boundaries of 
stabilizers of adjacent vertices, induced by edge stabilizers.
 Here the relation $\sim$ consists in glueing each peripheral circle of a Sierpinski carpet on a peripheral circle of another Sierpinski carpet. 
From this characterisation, we deduce that, the boundary of $\pi_1(X)$ is homeomorphic to the boundary of $\pi_1(X')$, where $X'$ is the graph of groups obtained from $X$ by replacing the vertex groups by some hyperbolic group whose boundary is a Sierpinski carpet, and by replacing the edge groups by the stabilizers of the peripheral circles. 
The result of such a construction is proven in \cite{KK} (Theorem 5) to be a 2-sphere, therefore  $\partial \Gamma$ is a 2-sphere. $\square \smallskip$

This proves Theorem \ref{theo;extension}, and even a little more since we learnt 
(Lemma \ref{lem;parab_in_C}) that the bounded parabolic points in a Sierpinski carpet cannot lie on the peripheral circles.

\subsection{Parabolic groups acting on Sierpinski carpets.}

We now prove Theorem  \ref{theo;intro2}.

{\it Proof: } Let $\Gamma$ be a relatively hyperbolic group whose boundary is a 
Sierpinski carpet, or a 2-sphere.
Let $H$ be a maximal parabolic subgroup of $\Gamma$. 
By Theorem \ref{theo;extension}, $H$ is a subgroup of 
finite index of a maximal parabolic subgroup of a 
relatively hyperbolic group whose boundary is a 2-sphere. 
Therefore we can assume without loss of generality that the boundary of 
$\Gamma$ is a 2-sphere $S^2$. 
By definition of bounded parabolic points, $H$ acts properly discontinuously 
co-compactly on the 
complement of a point in $S^2$, that is, on the plane 
$\mathbb{R}^2$. 
Since it is a proper action, the kernel is finite. The quotient is then a 2-dimensional good orbifold with infinite fundamental group. It is known (see \cite{Thu}, Chap.~XIII) that it is either Euclidean or hyperbolic. In both cases, the fundamental group is a virtual surface group. $\square$

 We can complete the result by the following.

\begin{prop}
  Let $(\Gamma,\mathcal{G})$ be a relatively hyperbolic group. Let $\mathcal{G}'$ be the family of subgroups of $\Gamma$ consisting of the elements of $\mathcal{G}$ that are not hyperbolic. Then, $(\Gamma, \mathcal{G}')$ is relatively hyperbolic.
\end{prop}
 
The proposition is a consequence of a characterisation of D.~Osin, Theorem 2.37 in \cite{O}. We can deduce:

\begin{coro}\label{coro;sierpinski_minimal}

Let $(\Gamma,\mathcal{G})$ be a  relatively hyperbolic group, such that $\mathcal{G}$ is minimal for this property, and  whose boundary is a Sierpinski carpet or a 2-sphere. Then every  element of $\mathcal{G}$ is virtually $\mathbb{Z} +\mathbb{Z}$.
\end{coro}

 There are examples of relatively hyperbolic groups with boundary homeomorphic
 to the Sierpinski carpet. Consider a geometrically finite hyperbolic 3-manifold whose convex core has at least 
 two boundary components, a cusp of rank $2$ and has no cusp of rank $1$. Let $\G$ be its fundamental group: $M=\mathbb{H}^3/\Gamma$. Let $B_1, \dots , B_n$ be the components of the boundary of the convex core. Their universal covers, and their translates by elements of $\G$, in $\mathbb{H}^3$, are disjoint hyperbolic planes that are not asymptotic to each other (such a thing would come from a cusp of rank $1$), and the fundamenta groups $\pi_1(B_i)$ are closed surface groups.
 The group $\Gamma$ is hyperbolic relative to its cusp subgroups and its boundary is its limit set $\Lambda \Gamma$ in $\partial \mathbb{H}^3$. It is then $S^2$ where a dense family of disjoint non-tangent open discs have been removed (each is stabilised by a certain conjugate of a group $\pi_1(B_i)$). See for example a similar study \cite{KK}.   Hence $\Lambda \Gamma$ is  homeomorphic to the
 Sierpinski carpet. Now consider $\G$ as hyperbolic relative to the cusp subgroups and the conjugates of the group $\pi_1
 (B_1)$. Its boundary is now a quotient of the previous  Sierpinski curve,
 where each peripheral circle fixed by a conjugate of $\pi_1(B_1)$ has been
 collapsed to a point. On the sphere $S^2 \simeq \partial \mathbb{H}_3$, this
 means contracting each disc stabilized by a conjugate of $\pi_1(B_1)$ to a
 point. But the result of such an operation is still a sphere, and, once
 removed the discs fixed by the conjugates of $\pi_1(B_j)$, $j\geq 2$, one still have a
 Sierpinski curve. The parabolic subgroups are then isomorphic either to
 $\mathbb{Z}^2$ (for the groups of the cusps) or to the surface group $\pi_1(B_1)$.


{\footnotesize

\thebibliography{99}

\bibitem{BGS}{\it W.~Ballmann, M.~Gromov, V.~Schroeder}, "Manifolds of nonpositive curvature", Progress in Mathematics, 61.
   Birkh\"auser (1985).


\bibitem{BM}{\it A.~Beardon, B.~Maskit}, ``Limit points of Kleinian groups, and finite sided fundamental polyhedra'', 
Acta Math. {\bf 132} (1974) 1-12.
 
\bibitem{Bparabolic}{\it B.H.~Bowditch}, "Discrete parabolic groups", 
J. Differential Geom. {\bf 38} (1993), no. 3, 559-583.

\bibitem{Brel}{\it B.H.~Bowditch}, ``Relatively hyperbolic groups'', preprint Southampton (1999).


\bibitem{BGF}{\it B.H.~Bowditch}, ``Boundaries of geometrically finite groups'', 
Math. Z. {\bf 230} (1999), 509-527.

\bibitem{BPeriph}{\it B.H.~Bowditch}, ``Peripheral splittings of groups'', Trans. Amer. Math. Soc. {\bf 353} (2001), 4057-4082. 

\bibitem{Dproc}{\it F.~Dahmani}, ``Classifying spaces and boundaries for relatively hyperbolic groups'', Proc. London Math. Soc. {\bf 86} (2003) 666-684.

\bibitem{Dcomb}{\it F.~Dahmani}, ``Combination of convergence groups'',  Geometry \& Topology {\bf 7} (2003) 933-963. 


\bibitem{Dsym}{\it F.~Dahmani, A.~Yaman} ``Symbolic dynamics for relatively hyperbolic groups'', preprint, (2002).

\bibitem{F}{\it B.~Farb}, ``Relatively Hyperbolic Groups'', Geom.Funct.Anal.{\bf 8} (1998), n$^o$ 5, 810-840.


\bibitem{Gab}{\it D.~Gabai}, ``Convergence groups are Fuchsian groups'', Ann. Math. {\bf 136} (1992) 447-510.


\bibitem{G}{\it M.~Gromov}, ``Hyperbolic groups'', Essays in group theory, MSRI, ed. S.Gersten, (1987), 75-263.

 

\bibitem{KK}{\it M.~Kapovich, B.~Kleiner}, ``Hyperbolic groups with low dimensional boundaries'', 
Ann. Scient. Ec. Norm. Sup. {\bf 33} 647-669.   


\bibitem{O}{\it D.~Osin}, ``Relatively hyperbolic groups: Intrinsic geometry, algebraic properties, and algorithmic problems'', Preprint 2003.

\bibitem{Szcz}{\it A.~Szczepa\'nski}, ``Relatively hyperbolic groups'' Michigan Math. J. {\bf 45} (1998) 611-618.

\bibitem{Thu}{\it W.~Thurston}, "The geometry and the topology of three-manifolds" 
Lecture notes, Princeton University, (1976-79). 

\bibitem{TukiaCrelle}{\it P.~Tukia}, ``Conical limit points and uniform convergence groups'',  J. reine angew. Math.  {\bf 501}  (1998), 71--98.

\bibitem{Tu}{\it P.~Tukia}, ``Generalizations of Fuchsian and Kleinian
groups'', First European Congress of Mathematics, Vol. II (Paris, 1992), 447-461,  Progress in Mathematics, Birkhauser (1994).



\bibitem{ref45deKK}{\it G.~Whyburn}, ``Topological characterization of the Sierpinski curve'', Fundamenta Math. {\bf 45} (1958) 320-324. 

}
 \end{document}